\documentclass[11pt]{article}

\usepackage{mathrsfs}
\usepackage{latexsym}
\usepackage{amsmath,amsthm}
\usepackage{amsfonts}
\usepackage{url}

\usepackage{amssymb}

\usepackage{mytheo}

\def\longto{\mathop{\longrightarrow}\limits}

\usepackage{bbm}
\newcommand{\chr}{\boldsymbol{\mathbbm{1}}} %
\newcommand{\pred}[1]{\chr_{\left\{ #1 \right\}}}

\newcommand{\E}{\mathbb{E}}
\renewcommand{\P}{\mathbb{P}}

\usepackage{geometry}
\makeatletter

\newcommand{\ben}{\begin{enumerate}}
\newcommand{\een}{\end{enumerate}}
\newcommand{\bit}{\begin{itemize}}
\newcommand{\eit}{\end{itemize}}

\renewcommand{\vec}[1]{\bs{\mathrm{#1}}}

\newcommand{\R}{\mathbb{R}}
\newcommand{\N}{\mathbb{N}}

\newcommand{\beq}{\begin{eqnarray*}}
\newcommand{\eeq}{\end{eqnarray*}}
\newcommand{\beqn}{\begin{eqnarray}}
\newcommand{\eeqn}{\end{eqnarray}}
\newcommand{\paren}[1]{\left( #1 \right)}
\newcommand{\sqprn}[1]{\left[ #1 \right]}
\newcommand{\tlprn}[1]{\left\{ #1 \right\}}
\newcommand{\set}[1]{\tlprn{#1}}

\newcommand{\bs}{\boldsymbol}

\newcommand{\hide}[1]{}
\newcommand{\oo}[1]{\frac{1}{#1}}

\def\eps{\varepsilon}

\newcommand{\evalat}[2]{\left.#1\right|_{#2}}

\title{
On the Concentration of the Missing Mass
}
\author{
Daniel Berend
and
Aryeh Kontorovich
}

\begin{document}
\maketitle
\begin{abstract}
A random variable is sampled from a discrete distribution. The missing mass
is the probability of the set of points not observed in the sample.
We sharpen and simplify McAllester and Ortiz's results (JMLR, 2003)
bounding the probability of large deviations of the missing mass.
Along the way, we refine and rigorously prove a fundamental inequality
of Kearns and Saul (UAI, 1998).

\end{abstract}

\newcommand{\lam}{\lambda}
\section{Introduction}
Hoeffding's classic inequality \cite{hoeffding} 
states that
If $X$ is a $[a,b]$-valued random variable with $\E X=0$ then
\beqn
\label{eq:hoeff}
\E e^{t X} \le e^{(b-a)^2t^2/8},
\qquad t\ge0
.
\eeqn
A standard proof of (\ref{eq:hoeff}) proceeds by writing
$x\in[a,b]$ as 
$x=pb+(1-p)a$, for $p=-a/(b-a)$,
and 
using convexity to obtain
\beqn
\label{eq:hoeff-proof}
\E e^{tX/(b-a)}\le 
(1-p)e^{-tp}+pe^{t(1-p)}
:=f(t)
\le e^{t^2/8}
,
\eeqn
where the last inequality follows
by noticing that $\log f(0)=\evalat{[\log f(t)]'}{t=0}=0$
and that $[\log f(t)]''\le1/4$.

Although (\ref{eq:hoeff}) is tight, it is a ``worst-case'' bound over
all distributions with the given support.
Refinements of (\ref{eq:hoeff})
include
the Bernstein and Bennett inequalities \cite{lugosi03},
which take the variance into account --- 
but these are also
too crude for some purposes.

In 1998, Kearns and Saul \cite{DBLP:conf/uai/KearnsS98}
put forth an exquisitely delicate inequality for (generalized) Bernoulli random variables,
which is sensitive
to the underlying distribution:
\beqn
\label{eq:k-s}
(1-p)e^{-tp}+pe^{t(1-p)}\le \exp\paren{\frac{1-2p}{4\log((1-p)/p)}t^2},
\qquad p\in[0,1],~t\in\R.
\eeqn
One easily verifies that (\ref{eq:k-s}) is superior to
(\ref{eq:hoeff}) ---
except for $p=1/2$, where the two coincide.
In fact, 
(\ref{eq:k-s}) is optimal in the sense that,
for every $p$, there is a 
$t$ for which equality is achieved.
The Kearns-Saul inequality
allows one to analyze various inference algorithms
in neural networks, and the influential 
paper \cite{DBLP:conf/uai/KearnsS98}
has inspired a fruitful line of research
\cite{
DBLP:journals/jair/BhattacharyyaK01,
DBLP:journals/jmlr/McAllesterO03,
DBLP:conf/nips/NgJ99,
DBLP:conf/nips/NguyenJ03}.

One specific application of the Kearns-Saul inequality involves the concentration of
the {missing mass}. Let $\vec p=(p_1,p_2,\ldots)$ be a distribution over $\N$
and suppose that $X_1,X_2,\ldots,X_n$ are sampled iid according to $\vec p$. Define the
indicator variable $\xi_j$ to be $0$ if $j$ occurs in the sample and $1$ otherwise:
\beq
\xi_j=\pred{j\notin\set{X_1,\ldots,X_n}},
\qquad j\in\N.
\eeq
The {\em missing mass} 
is 
the random variable
\beqn
\label{eq:Undef}
U_n = \sum_{j\in\N} p_j\xi_j.
\eeqn
McAllester and Schapire \cite{DBLP:conf/colt/McAllesterS00}
first established subgaussian concentration for the missing mass
via a somewhat intricate argument. Later, McAllester and Ortiz
\cite{DBLP:journals/jmlr/McAllesterO03}
showed how the standard inequalities of
Hoeffding, Angluin-Valiant, Bernstein and Bennett
are inadequate for obtaining 
exponential bounds of the correct
order in $n$,
and developed a thermodynamic
approach for systematically handling this problem\footnote{
The latter has, in turn, inspired a general thermodynamic approach to concentration
\cite{pre06038702}.
}.

We were led to the Kearns-Saul inequality (\ref{eq:k-s})
in an attempt to understand and simplify the missing mass concentration results
of 
McAllester and Ortiz 
\cite{DBLP:journals/jmlr/McAllesterO03},
some of which rely on (\ref{eq:k-s}). However, we were unable to complete the proof of (\ref{eq:k-s}) sketched in \cite{DBLP:conf/uai/KearnsS98},
and a literature search likewise came up empty. 
The proof we give here follows an alternate path, and may be of independent interest.
As an application,
we 
simplify and sharpen some of the missing mass concentration results given in \cite{DBLP:conf/colt/McAllesterS00,DBLP:journals/jmlr/McAllesterO03}.

\section{Main results}
In \cite[Lemma 1]{DBLP:conf/uai/KearnsS98}, Kearns and Saul define the function
\beqn
g(t) = \oo{t^2} \log\sqprn{
(1-p)e^{-tp}+pe^{t(1-p)}},
\qquad t\in\R
.
\eeqn
A natural attempt to find the maximum of $g$ leads one
to
the transcendental equation
\beq
g'(t)=
\frac{
(e^t-1)(1-p)pt-2(1+(e^t-1)p)\log[1+(e^t-1)pe^{-pt}]
}{
(1+(e^t-1)p)t^3
}=0.
\eeq
In an inspired tour de force, Kearns and Saul were able to find that $g'(t^*)=0$
for
\beq
t^* = 2\log\frac{1-p}{p}.
\eeq
This observation naturally suggests (i) arguing that $t^*$ is the unique zero of $g'$
and (ii) supplying (perhaps via second-order information) an argument for $t^*$ being a local maximum.
In fact, all evidence points to $g'(t)$ having the following properties:
\bit
\item[(*)] $g'>0$ on $(-\infty,t^*)$,
\item[(**)] $g'=0$ at $t=t^*$,
\item[(***)] $g'<0$ on $(t^*,\infty)$.
\eit
Unfortunately, besides straightforwardly verifying (**), we were not able to formally establish
(*) or (***) --- and we leave this as an intriguing open problem.
Instead, in Theorem \ref{thm:main} we prove the Kearns-Saul inequality (\ref{eq:k-s})
via a rather different approach. Moreover, for $p\ge1/2$ and $t\ge0$,
the right-hand side of (\ref{eq:k-s}) may be improved to
$\exp[p(1-p)t^2/2]$. This refinement, proved in Lemma \ref{lem:k-s-refine},
may be of independent interest.

As an application, 
we recover the upper tail estimate on the missing mass in
\cite[Theorem 16]{DBLP:journals/jmlr/McAllesterO03}:
\bethn
\label{thm:upper-tail}
\beq
\P(U_n>\E U_n+\eps)\le e^{-n\eps^2}.
\eeq
\enthn
We also obtain the following 
lower tail estimate:
\bethn
\label{thm:lower-tail}
\beq
\P(U_n<\E U_n-\eps)\le e^{-C_0n\eps^2/4},
\eeq
where
\beq
C_0=\inf_{0<x<1/2}
~
\frac{2}{x(1-x)\log(1/x)}
\approx7.6821.
\eeq
\enthn
Since $C_0/4\approx1.92$, 
Theorem \ref{thm:lower-tail}
sharpens the estimate in
\cite[Theorem 10]{DBLP:journals/jmlr/McAllesterO03},
where the constant in the exponent was $e/2\approx 1.36$.
Our bounds are arguably simpler than those in 
\cite{DBLP:journals/jmlr/McAllesterO03} as they bypass the thermodynamic approach.

\section{Proofs}
The following well-known estimate is an immediate consequence of (\ref{eq:hoeff-proof}):
\belen
\label{lem:cosh}
\beq
\oo2e^{-t}+\oo2e^{t}=\cosh t \le e^{t^2/2},
\qquad t\in\R
.
\eeq
\enlen

We proceed with a proof of the Kearns-Saul inequality.
\bethn
\label{thm:main}
For all $p\in[0,1]$ and $t\in\R$,
\beqn
\label{eq:main}
(1-p)e^{-tp}+pe^{t(1-p)}\le \exp\paren{\frac{1-2p}{4\log((1-p)/p)}t^2}.
\eeqn
\enthn
\bepf
The cases $p=0,1$ are trivial.
Since
\beq
\lim_{p\to1/2} \frac{1-2p}{\log((1-p)/p)}=1/2,
\eeq
for $p=1/2$ the claim 
follows from Lemma \ref{lem:cosh}.
For $p\neq1/2$, 
we multiply both sides of 
(\ref{eq:main})
by $e^{tp}$, take logarithms,
and put $t=2s\log((1-p)/p)$ to obtain the equivalent
inequality
\beqn
\label{eq:t->s}
s\paren{ s+2p(1-s)}\log((1-p)/p)
-\log\paren{ 1-p+p( (1-p)/p)^{2s}}
\ge0.
\eeqn
For $s\in\R$, denote the left-hand side of (\ref{eq:t->s}) by $h_s(p)$.
A routine calculation yields
\beqn
\label{eq:h_s(1/2)}
h_s(1/2)=h_s'(1/2)=0
\eeqn
and
\beq
h_s''(p)=
\paren{
\frac{ (\mu-1)p^2-s+p(1-\mu+s+\mu s)}
{ p(1-p)(1+(\mu-1)p) }
}^2,
\eeq
where $\mu=((1-p)/p)^{2s}$.

As $h_s''\ge0$, 
we have that
$h_s$ is convex, and 
from (\ref{eq:h_s(1/2)})
it follows that $h_s(p)\ge0$ for all $s,p$.
\enpf

We will also need a refinement of (\ref{eq:k-s}):
\belen
\label{lem:k-s-refine}
For $p\in[1/2,1]$ and $t\ge0$,
\beqn
\label{eq:k-s-refine}
\oo{t^2}\log\sqprn{
(1-p)e^{-tp}+pe^{t(1-p)}}
\le
\frac{p(1-p)}{2}.
\eeqn
\enlen
{\bf Remark:} Since the right-hand side of (\ref{eq:main}) majorizes
the right-hand side of (\ref{eq:k-s-refine}) uniformly over $[1/2,1]$, the latter estimate
is tighter.
\bepf
The claim is equivalent to
\beq
L(p):=(1-p)+pe^t\le \exp(pt+p(1-p)t^2/2)=:R(p),
\qquad t\ge0.
\eeq
For $t\ge4$, we have
\beq
\exp(pt+p(1-p)t^2/2) &\ge&
\exp(pt+p(1-p)(4t)/2)\\ 
&=&
\exp(pt+2p(1-p)t) \\
&\ge&
\exp(pt+(1-p)t)
= e^t\ge L(p).
\eeq
For $0\le t <4$,
\beq
R''(p)-L''(p) = 
\oo4 \exp(pt(2+t-pt)/2)(2p-1)t^3( (2p-1)t-4),
\qquad p\in[1/2,1],
\eeq
which is obviously non-positive.
Now the inequality clearly holds at $p=1$ (as equality),
and the $p=1/2$ case is 
implied by Lemma \ref{lem:cosh}.
The claim now follows by convexity.

\enpf

Our numerical constants are defined in the following lemma,
whose elementary proof is omitted:
\belen
\label{lem:x0}
Define the function
\beq
f(x) = 
{x(1-x)\log(1/x)},
\qquad
x\in(0,1/2).
\eeq
Then 
$x_0\approx0.2356$ 
is the unique solution of
$f(x)'=0$ on $(0,1/2)$.
Furthermore, 
\beqn
\label{eq:c0}
C_0:=
\inf_{0<x<1/2}2/f(x)
=2/f(x_0)
\approx7.6821.
\eeqn
\enlen

The main technical step towards obtaining our 
missing mass deviation estimates is
the following lemma.
\belen
\label{lem:lamp}
Let $n\ge1$, $\lam\ge0$, $p\in[0,1]$, 
and put $q=(1-p)^n$.
\hide{
Let
$x_0\approx0.2356$ 
be the (unique) solution of 
the equation
\beqn
\label{eq:x0}
\log x=-\frac{1-x}{1-2x},
\qquad
x\in[0,1/2]
\eeqn
and put $C_0=2/[x_0(1-x_0)\log(1/x_0)]\approx7.6821$.
}
Then:
\bit
\item[(a)]
$$ 
q e^{\lam(p-pq)}+(1-q)e^{-\lam p q}
 \le \exp( p\lam^2/4n).
$$
\item[(b)]
$$ 
q e^{\lam(pq-p)}+(1-q)e^{\lam p q}
 \le \exp( p\lam^2/C_0n).
$$
\eit
\enlen
\bepf
\bit
\item[(a)]
We invoke Theorem \ref{thm:main}
with
$p=q$ and
$t=\lam p$ 
to obtain
\beq
q e^{\lam(p-pq)}+(1-q)e^{-\lam p q}
\le
\exp[ (1-2q)\lam^2p^2/4\log[(1-q)/q]].
\eeq
Thus it suffices to show that
\beq
(1-2q)\lam^2p^2/4\log[(1-q)/q] \le
p\lam^2/4n,
\eeq
or equivalently,
\beq
(1-2q)p/\log[(1-q)/q] \le
\log(1-p)/\log q,
\qquad p,q\in[0,1]
.
\eeq
Collecting the $p$ and $q$ terms on opposite sides, 
it remains to prove that
\beq
L(q):=\frac{(1-2q)\log(1/q)}{\log[(1-q)/q]}
\le
\frac{\log(1/(1-p))}{p}=:R(p),
\qquad 0<p,q<1.
\eeq
We claim that 
$L\le1\le R$.
The 
second inequality
is obvious from the Taylor expansion, since
\beqn
\label{eq:p-fac>1}
\frac{\log(1/(1-p))}{p}
=
1+p/2+p^2/3+p^3/4+\ldots.
\eeqn
To prove 
that $L\le1$,
we note first that $L(q)\ge L(1-q)$ for $q\in(0,1/2)$.
Hence, it suffices to consider $q\in(0,1/2)$.
To this end,
it suffices to show that
the function
\beq
f(q) = \log[(1-q)/q]-(1-2q)\log(1/q)
\eeq
is positive on $(0,1/2)$. 
Since $\lim f(q) \longto_{q\to0}0=f(1/2)$ and
\beq
f''(q) = \frac{ -2+3q-2q^2}{(1-q)^2q}\le0,
\eeq
it follows that
$f\ge0$ on $[0,1/2]$.

\item[(b)]
The inequality is equivalent to 
\beq
L(\lam):=
\oo{\lam^2 p^2}
\log\sqprn{ qe^{-\lam p(1-q)}+(1-q)e^{\lam p q}}
\le \oo{\lam^2 p^2}\frac{\lam^2 p}{C_0\log q/\log(1-p)}
=:R
,
\eeq
where $L$ is obtained
from the left-hand side of (\ref{eq:main})
after replacing $p$ by $1-q$ and $t$ by $\lam p$.
We analyze the cases 
$q<1/2$ and 
$q>1/2$ separately (as above, the case where $q=1/2$ is trivial).
For $q>1/2$,
put $\lam^*=\frac{2}{p}\log\frac{q}{1-q}>0$
and invoke 
Theorem \ref{thm:main}
to conclude that
$\sup_{\lam\ge0}L(\lam)\le L(\lam^*)$.
Hence,
it remains to prove that 
$L(\lam^*)\le R$, or equivalently,
\beq
\frac{(2q-1)}{4\log(q/(1-q))} \le \oo{\lam^2 p^2}\frac{\lam^2 p}{C_0\log q/\log(1-p)}
.
\eeq
After simplifying, 
this amounts to showing
that 
\beq
4 \frac{\log(1/(1-p))}{p} \frac{\log(q/(1-q))}{(2q-1)\log(1/q)}
\ge C_0.
\eeq

As 
in (\ref{eq:p-fac>1}),
the factor 
${\log[1/(1-p)]}/{p}$
is bounded below by $1$. 
We claim that the 
factor $\frac{\log(q/(1-q))}{(2q-1)\log(1/q)}$,
increases for $q\in[1/2,1]$.
Indeed, this is obvious for $1/\log(1/q)$, and the expansion
about $q=1/2$
\beq
\frac{\log(q/(1-q))}{(2q-1)}
= \sum_{n=0}^\infty \frac{2^{2n+1}}{2n+1}\paren{q-\oo2}^{2n}
\eeq
shows that the same holds for $\frac{\log(q/(1-q))}{2q-1}$.
In particular,
\beq
4 \frac{\log(1/(1-p))}{p} \frac{\log(q/(1-q)}{(2q-1)\log(1/q)
}
&\ge&
4\cdot1\cdot \lim_{q\to1/2} \frac{\log(q/(1-q)}{(2q-1)\log(1/q)
}\\ 
&=& 8/\log2\approx 11.542>C_0.
\eeq

When $q<1/2$,
we invoke Lemma \ref{lem:k-s-refine}
together with the observation that
\beq
\lim_{\lam\to0_+}L(\lam)= \frac{q(1-q)}{2}
\eeq
to conclude that
$\sup_{\lam\ge0}L(\lam)\le L(0)$.
Hence,
it remains to show that 
\beq
\frac{\log(1/(1-p))}{p}
\frac{2}{q(1-q)\log(1/q)}
\ge C_0.
\eeq
As in (\ref{eq:p-fac>1}),
${\log[1/(1-p)]}/{p}\ge1$
and
the claim follows
by Lemma \ref{lem:x0}.
\eit
\enpf

Our proof of Theorems \ref{thm:upper-tail} and \ref{thm:lower-tail}
is facilitated by the following observation,
also made in \cite{DBLP:journals/jmlr/McAllesterO03}.
Although the random variables $\xi_j$ whose weighted sum comprises the
missing mass (\ref{eq:Undef}) are not independent,
they are {\em negatively associated} \cite{Dubhashi:1998:BBS:299633.299634}.
A basic fact about negative association is that it is
``at least as good as independence'' as far as exponential concentration is
concerned 
\cite[Lemmas 5-8]{DBLP:journals/jmlr/McAllesterO03}:
\belen
\label{lem:xi'}
Let $\xi_j'$ be independent random variables,
where $\xi_j'$ is distributed identically to $\xi_j$
for all $j\in\N$. Define also the ``independent analogue''
of $U_n$:
\beq
U_n'=\sum_{j\in\N} p_j\xi_j'.
\eeq
Then
for all $n\in\N$ and $\eps>0$,
\bit
\item[(a)]
\beq
\P(U_n\ge\E U_n+\eps)&\le& \P(U_n'\ge\E U_n'+\eps),
\eeq
\item[(b)]
\beq
\P(U_n\le\E U_n-\eps)&\le& \P(U_n'\le\E U_n'-\eps).
\eeq
\eit
\enlen
\bepf[Proof of Theorems \ref{thm:upper-tail} and \ref{thm:lower-tail}]
Observe that 
the random variables $\xi_j'$
defined in Lemma \ref{lem:xi'}
have
a Bernoulli distribution
with $\P(\xi_j'=1)=q_j=(1-p_j)^n$ and
put
$X_j=\xi_j-\E\xi_j$.
Using standard exponential bounding with Markov's inequality,
\beq
\P(U_n\ge\E U_n+\eps)&\le& \P(U_n'\ge\E U_n'+\eps) \\
&=& P\sqprn{\exp\paren{\lam\sum_{j\in\N}X_j}\ge e^{\lam\eps}},
\qquad \lam\ge0\\
&\le& e^{-\lam\eps}\prod_{j\in\N}\E e^{\lam X_j}\\
&=& e^{-\lam\eps}\prod_{j\in\N}
\paren{q_j e^{\lam(p_j-p_jq_j)}+(1-q_j)e^{-\lam p_j q_j}}\\
&\le&
e^{-\lam\eps}\prod_{j\in\N}
\exp( p_j\lam^2/4n)\\
&=&\exp(\lam^2/4n-\lam\eps),
\eeq
where the last inequality invoked Lemma \ref{lem:lamp}(a).
Choosing $\lam=2n\eps$ yields Theorem
\ref{thm:upper-tail}.

The proof of the 
Theorem
\ref{thm:lower-tail}
is almost identical,
except that 
$X_j$ is replaced by $-X_j$ and
Lemma \ref{lem:lamp}(b) is invoked
instead of Lemma \ref{lem:lamp}(a).
\enpf

\bibliographystyle{plain}
\bibliography{../mybib}

\end{document}